\documentclass[12pt,leqno]{article}
\usepackage[pdftex]{graphicx}
\usepackage{amsfonts, amsmath}
\usepackage{times,wrapfig,verbatim}
\usepackage{color}

\newcounter{conjecture}
\newcounter{remark}
\newcounter{corollary}
\newenvironment{conjecture}{\medskip{\bf Conjecture.\ }\em}{\rm}

\newenvironment{corollary}{\medskip{\bf Corollary \thecorollary.}
\addtocounter{corollary}{1}\em}{\rm}
\newtheorem{theorem}{Theorem}
\newtheorem{lemma}{Lemma}
\newtheorem{proposition}{Proposition}

\newcommand{\dd}{\delta}

\newcommand{\aaa}{\alpha}

\newcommand{\lll}{\label}
\newcommand {\rrr}[1]{(\ref{#1})}

\def \be{\begin{equation}}
\def \ee{\end{equation}}
\def \bt{\begin{theorem}}
\def \et{\end{theorem}}
\def \bc{\begin{corollary}}
\def \ec{\end{corollary}}
\def \bea{\begin{eqnarray}}
\def \eea{\end{eqnarray}}
\def \bas{\begin{eqnarray*}}
\def \eas{\end{eqnarray*}}

\def \aa{\alpha}
\def \bb{\beta}
\def \ga{\gamma}
\def \Ga{\Gamma}
\def \GGG{\Gamma}

\def \la{\lambda}

\def \th{\theta}

\def \vski{\vspace{12pt}}

\def \({\left(}
\def \){\right)}

\def \bc{\begin{center} }
\def \ec{\end{center} }
\def \bs{\begin{slide} }
\def \es{\end{slide} }

\def\square{{\vcenter{\vbox{\hrule height.3pt
         \hbox{\vrule width.3pt height5pt \kern5pt
            \vrule width.3pt}
         \hrule height.3pt}}}}
\def\qed{{\hfill $\square$ \bigskip}}

\newcounter{cccases}
\begin{document}

\title{On the Cheeger constant for distance-regular graphs.}


 \author{Jack H. Koolen$^{\,\rm 1}$,~Greg Markowsky$^{\,\rm 2}$ and Zhi Qiao$^{\,\rm 3}$\\
{\small {\tt koolen@ustc.edu.cn} ~~ {\tt gmarkowsky@gmail.com} ~~ {\tt zhiqiao@sicnu.edu.cn}}\\
{\footnotesize{$^{\rm 1}$Wen-Tsun Wu Key Laboratory of the CAS and School of Mathematical Sciences, USTC, Hefei, China}}\\
{\footnotesize{$^{\rm 2}$Department of Mathematics,  Monash University, Australia}}\\
{\footnotesize{$^{\rm 3} $College of Mathematics and Software Science, Sichuan Normal University, China}}}

\maketitle

\begin{abstract}
The Cheeger constant of a graph is the smallest possible ratio between the size of a subgraph and the size of its boundary. It is well known that this constant must be at least $\frac{\la_1}{2}$, where $\la_1$ is the smallest positive eigenvalue of the Laplacian matrix. The subject of this paper is a conjecture of the authors that for distance-regular graphs the Cheeger constant is at most $\la_1$. In particular, we prove the conjecture for the known infinite families of distance-regular graphs, distance-regular graphs of diameter 2 (the strongly regular graphs), several classes of distance-regular graphs with diameter 3, and most distance-regular graphs with small valency.
\end{abstract}

\section{Introduction}

The {\it Cheeger constant} $h_G$ of a graph $G$ is a prominent measure of the connectivity of $G$, and is defined as

\begin{equation} \label{}
h_G = \inf \Big\{ \frac{E[S,S^c]}{\mbox{vol}(S)}| S \subset V(G) \mbox{ with } |S| \leq \frac{|V(G)|}{2}\Big\},
\end{equation}

where $V(G)$ is the vertex set of $G$, $\mbox{vol}(S)$ is the sum of the valencies of the vertices in $S$, $S^c$ is the complement of $S$ in $V(G)$, $|S|$ is the number of vertices in $S$, and for any sets $A,B$ we use $E[A,B]$ to denote the number of edges in $G$ which connect a point in $A$ with a point in $B$; to make this last definition more precise, let $E(G) \subseteq V(G) \times V(G)$ denote the edge set of $G$, and then $E[A,B] = |\{ (x,y) \in E(G) \mbox{ with } x \in A, y \in B\}|$; note that with this definition each edge $(x,y)$ with $x,y \in A \cap B$ will in essence be counted twice, and $\mbox{vol}(S) = E[S,S] + E[S,S^c]$. There is an interesting connection between $h_G$ and spectral graph theory, given by the following general result.

\begin{theorem} \label{mass}
Let $\la_1$ be the smallest positive eigenvalue of the Laplacian matrix of $G$. Then

\begin{equation} \label{jc}
\frac{\la_1}{2} \leq h_G \leq \sqrt{\la_1(2-\la_1)}.
\end{equation}
\end{theorem}

See \cite{chung} for a proof of this statement. 
In this paper, we are interested in studying the Cheeger constant for the family of {\it distance-regular graphs} (see Section \ref{drgprop} for definitions). In particular, we make and provide evidence for the following conjecture.

\begin{conjecture} \label{}
Suppose $G$ is a distance-regular graph, and $\la_1$ is the smallest positive eigenvalue of the Laplacian matrix of $G$. Then

\begin{equation} \label{}
\frac{\la_1}{2} \leq h_G \leq \la_1.
\end{equation}
\end{conjecture}

We prove the conjecture for strongly regular graphs (distance-regular graphs with diameter 2), and for a number of families of distance-regular graphs and special cases. We should concede that there are several graphs for which we have been unable to verify the conjecture (for instance the flag graph of $GH(2,2)$ and incidence graph of $GH(3,3)$, see the end of Section \ref{smallval}), and for which the conjecture may fail. However, in that case we still conjecture that there are only finitely many distance-regular graphs for which $h_G > \la_1$, and hope that they can be classified.

\vski


We should mention that bounding the Cheeger constant on distance-regular graphs has already been considered, in \cite{krap}. However, we were unable to follow the methods given there and found a counterexample to their initial claim, which the authors acknowledged in the subsequent corrigendum. In any event, the methods we use are entirely different from theirs and for most graphs the inequality we obtain is stronger than the one they claim. The topic has also been mentioned in \cite{sivasiva}, although different types of questions than ours were addressed there.

\vski

In the next section, we provide the required definitions and notation. Section \ref{infam} proves the conjecture for the principal known infinite families of distance-regular graphs, Section \ref{srg} proves it for strongly regular graphs, Section \ref{diam3} proves it for several subclasses of distance-regular graphs of diameter 3, and Section \ref{smallval} proves it for most distance-regular graphs of small valency.

\section{Distance-regular graphs} \lll{drgprop}

All the graphs considered in this paper are finite, undirected and
simple (for unexplained terminology and more details, see for example \cite{drgraphs}). Let $G$ be a
connected graph and let $V(G)$ be the vertex set of $G$. The {\it distance} $d(x,y)$ between
any two vertices $x,y$ of $G$
is the length of a shortest path between $x$ and $y$ in $G$. The {\it diameter} of $G$ is the maximal distance
occurring in $G$ and we will denote this by $D = D(G)$.  
For a vertex $x \in V(G)$, define $\Ga_i(x)$ to be the set of
vertices which are at distance $i$ from $x~(0\le i\le
D)$, and when the choice of $x$ is unimportant we will simply write $\Ga_i$. In addition, define $\Ga_{-1}(x)=\Ga_{D+1}(x)
:= \emptyset$. We write $x\sim_{G} y$ or simply $x\sim y$ if two vertices $x$ and $y$ are adjacent in $G$. A connected graph $G$ with diameter $D$ is called
{\em distance-regular} if there are integers $b_i,c_i$ $(0 \le i
\leq D)$ such that for any two vertices $x,y \in V(G)$ with $d(x,y)=i$, there are precisely $c_i$
neighbors of $y$ in
$\Ga_{i-1}(x)$ and $b_i$ neighbors of $y$ in $\Ga_{i+1}(x)$
(cf. \cite[p.126]{drgraphs}). In particular, a distance-regular graph $G$ is regular with valency
$k := b_0$ and we define $a_i:=k-b_i-c_i$ for notational convenience. 
The numbers $a_i$, $b_{i}$ and $c_i~(0\leq i\leq D)$ are called the {\em
intersection numbers} of $G$, and the sequence $\{b_0,b_1, \ldots, b_{D-1};c_1,c_2, \ldots , c_D\}$ is the {\it intersection array} of $G$. Note that always $b_D=c_0=a_0=0$, $b_0 = k$ and $c_1=1$.
The intersection numbers of a distance-regular graph $G$ with diameter $D$ and valency $k$ satisfy
(cf. \cite[Proposition 4.1.6]{drgraphs})\\

(i) $k=b_0> b_1\geq \cdots \geq b_{D-1}$;\\
(ii) $1=c_1\leq c_2\leq \cdots \leq c_{D}$;\\
(iii) $b_i\ge c_j$ \mbox{ if }$i+j\le D$.\\

Moreover, if we fix a vertex $x$ of $G$, then $|\Ga_i(x)|$ does not depend on the
choice of $x$ as $c_{i+1} |\Ga_{i+1}(x)| =
b_i |\Ga_i(x)|$ holds for $i =1, 2, \ldots, D-1$.

\vski

For a distance-regular graph $G$ of diameter $D$, we will write $k=\th_0 > \th_1 > \ldots > \th_D$ to describe the eigenvalues of the adjacency matrix $A$ of $G$, and refer to $\th_0, \ldots , \th_D$ as simply the {\it eigenvalues of $G$}. The {\it Laplacian matrix} (sometimes referred to as the {\it normalized Laplacian}) $L = I - \frac{1}{k}A$ will therefore have eigenvalues $0=\la_0 < \la_1 < \ldots < \la_D$, where the relationship $\la_i = \frac{k-\th_i}{k}$ holds. We will refer to $\la_0, \ldots, \la_D$ as the {\it Laplacian eigenvalues of $G$}.

\section{Infinite families} \label{infam}

In this section we prove the conjecture for the principal known infinite families of distance-regular graphs with unbounded diameter. These families include many of the most well-known distance-regular graphs. In every case we will find an induced subgraph $G'$ which satisfies the conditions of the following lemma, which is simply a restatement of our conjecture in a form that is somewhat easier to check.

\begin{lemma} \label{iran}
If a regular graph $G$ on $v$ vertices with valency $k$ admits an induced subgraph $G'$ on $v' \leq \frac{v}{2}$ vertices with average valency $k' = \frac{E[G',G']}{|G'|}$ where $k' \geq \th_1(G)$, then $h_G \leq \la_1$.
\end{lemma}

{\bf Proof:} In this case, $h_G \leq \frac{E[G',(G')^c]}{k|G'|} = \frac{k|G'|-E[G',G']}{k|G'|} = \frac{(k-k')|G'|}{k|G'|} \leq \frac{k-\th_1}{k} = \la_1$. \qed


The infinite families are as follows.

\begin{enumerate}
	\item $G$ is a Johnson graph $J(n,e)$ with $n\geq 2e \geq 2$ (\cite[p.255]{drgraphs}). The vertices of $G$ can be realized as $e$-subsets of the set $N=\{1, \ldots, n\}$. Then $\theta_1=(e-1)(n-e-1)-1$.
	If $G'$ is the induced subgraph on the set of all vertices that contain the element $1 \in N$, then  $G'$ is isomorphic to $J(n-1,e-1)$ with valency $k'=(e-1)(n-e) \geq \theta_1$. Also, $v=|G|=\binom{n}{e}=\frac{n}{e}\binom{n-1}{e-1}=\frac{n}{e}|G'|\geq 2|G'|$.
	Thus, Lemma \ref{iran} applies. 
	
	\item $G$ is a Hamming graph $H(d,q)$ (\cite[p.261]{drgraphs}).The vertices of $G$ can be realized as the elements of $N^d$, where $N=\{1, \ldots, q\}$, with two vertices being adjacent when they differ in exactly one component. Here $\theta_1=q(d-1)-d$. If $G'$ is the induced subgraph on the set of all vertices with first component equal to $1 \in N$, then $G'$ is isomorphic to $H(d-1,q)$ with valency $k'=(d-1)(q-1) \geq \theta_1$. Also, $v=|G|=q|G'|$.
	Thus, Lemma \ref{iran} applies. 

	\item $G$ is a Doob graph $G(d_1,d_2)$, which is the Cartesian product of $H(d_2,4)$ with $d_1$ copies of the Shrikhande graph, which has the same intersection numbers (and thus same spectrum) as $H(2d_1+d_2,4)$ (\cite[p.262]{drgraphs}). We therefore have $\theta_1=4(2d_1+d_2-1)-(2d_1+d_2) =6d_1+3d_2-4$. If $d_2>0$ then we take $G'$ to be the Cartesian product of the subgraph used for the Hamming graph before, which is isomorphic to $H(d_2-1,4)$, with $d_1$ copies of the Shrikhande graph. $G'$ has the same valency as $H(2d_1+d_2-1,4)$, and this is $3(2d_1+d_2-1) > \theta_1$. Also $|G| = 4|G'|$, so Lemma \ref{iran} applies. 
	Suppose $d_2=0, d_1 > 0$, so that $G$ is the Cartesian product of $d_1$ copies of the Shrikhande graph. The Shrikhande graph has valency $6$ and is locally a hexagon, and thus contains a 6-wheel (which is a hexagon with an additional point added, where the additional point is connected to every point in the hexagon). We can therefore take $G'$ to be the Cartesian product of a 6-wheel with $d_1-1$ copies of the Shrikhande graph. Then the vertices in $G'$ have valencies either $6(d_1 - 1) +6$ or $6(d_1 - 1) +3$, both of which are larger than $\th_1 = 6d_1 - 4$. Thus the average valency of $G'$ is greater than $\th_1$, and furthermore $|G'| = \frac{7}{16}|G| \leq \frac{1}{2}|G|$, so Lemma \ref{iran} applies. 

	\item $G$ is a halved $n$-cube, which is the bipartite half of the hypercube $H(n,2)$ (\cite[p.264]{drgraphs}). Here $\theta_1=\frac{1}{2}(n-2)^2-\frac{1}{2}n$.
	The vertices of $G$ can be realized as the set of all binary strings of length $n$ containing an even number of $1$'s, where two vertices are adjacent if they are of Hamming distance 2. Take the subgraph $G'$ to be induced on the set of all such strings with first digit $0$. $G'$ is therefore isomorphic to the halved $(n-1)$-cube, and therefore has valency $k'=\frac{1}{2}(n-1)(n-2)$.
	Then $k'\geq \theta_1$, and $|G| = 2 |G'|$, so Lemma \ref{iran} applies.

    \item $G$ is a folded $n$-cube, so $G$ is the graph $H(n-1,2)$ with a perfect matching introduced between antipodal vertices (\cite[p.264]{drgraphs}). Here $\th_1 = n-4$. We can again consider the vertices of $G$ to be the set of all binary strings of length $n-1$, and let the subgraph $G'$ be induced on the set of all such strings with first digit $0$; $G'$ is then isomorphic to an $(n-2)$-cube, so that $|G|=2|G'|$ and the valency of $G'$ is $n-2$. Thus, Lemma \ref{iran} applies.

    \item $G$ is a folded halved 2$n$-cube, which can be realized as the halved $2n$-cube discussed earlier with antipodal points identified(\cite[p.265]{drgraphs}). Here $\th_1 = 2(n-2)^2-n$. The antipodal identification means that one of the two strings associated to each vertex ends with a 0, and the vertices of $G$ can therefore be realized as the set of all binary strings of length $2n-1$ containing an even number of $1$'s. We may let $G'$ be the subgraph induced on the set of all such strings with first two digits $00$, and then it may be checked that $G'$ is isomorphic to a halved $(2n-3)$-cube and the valency of $G'$ is $k'=\frac{(2n-3)(2n-4)}{2}= 2n^2 - 7n + 6 \geq 2n^2 - 9n +8 = \th_1$. Furthermore $|G| \geq 2|G'|$, so Lemma \ref{iran} applies.  

    \item $G$ is the odd graph $O_k$, so the vertices of $G$ can be realized as all binary strings of length $2k-1$ containing exactly $(k-1)$ 1's, where two vertices are adjacent if the digit-wise product of their corresponding strings is 0 (that is, there are no digits where both strings are equal to 1). Here $k$ is the valency of the graph, and $\th_1 = k-2$ (\cite[p.260]{drgraphs}). Suppose $k \geq 3$, as the other cases are trivial. Let $A$ be the set of all strings in $G$ that begin with $1100$, and let $B$ be the set of all strings in $G$ that begin with $0011$. Let $G'$ be the subgraph induced on $A \cup B$. We see that $a \in A$ is adjacent to $b$ in $G'$ iff $b \in B$ and the digit-wise product of $\tilde a$ and $\tilde b$ is 0, where $\tilde a, \tilde b$ are the strings $a,b$ with the first four digits removed. It follows that $G'$ is isomorphic to the bipartite double cover of the odd graph $O_{k-2}$ (i.e. the doubled odd graph), which we will denote $(dO)_{k-2}$. The valency of $(dO)_{k-2}$ is the same as $O_{k-2}$, namely $k-2 = \th_1$. Also, $|G'| = 2|O_{k-2}| = 2 {2k-5 \choose k-3} = {2k-4 \choose k-2}$, so $\frac{|G|}{|G'|} = \frac{{2k-1 \choose k-1}}{{2k-4 \choose k-2}} = \frac{(2k-1)(2k-2)(2k-3)}{k(k-1)^2}>2$. Thus Lemma \ref{iran} applies.

	\item $G$ is a doubled odd graph $(dO)_{m+1}$, so we can realize the vertices of $G$ as pairs $\{X,v\}$, where $X \in \{A,B\}$ and $v$ is a binary string of length $2m+1$ containing exactly $m$ 1's. Two vertices $\{X_1,v_1\}$ and $\{X_2,v_2\}$ are adjacent iff $X_1 \neq X_2$ and the digit-wise product of $v_1$ and $v_2$ is 0 (\cite[p.260]{drgraphs}). Here $k=m+1$, and since the smallest eigenvalue of the odd graph $O_{m+1}$ is $-m$ (\cite[p.260]{drgraphs}), it follows from \cite[Thm. 1.11.1 (v)]{drgraphs} that $\th_1 = m$. Let $A'$ be the set of all vertices associated with strings in $A$ whose corresponding string starts with "01", and $B'$ be the set of all vertices associated with strings in $B$ whose corresponding string starts with "10". Let $G'$ be the subgraph induced on $A' \cup B'$. If $a' \in A'$, then in $G'$ the vertex $a'$ is adjacent to all vertices $b'$ in $B'$ such that the digit-wise product of $\tilde a'$ and $\tilde b'$ is 0, where $\tilde a'$ and $\tilde b'$ are the strings corresponding to $a'$ and $b'$ with the first two elements removed. The reverse statement holds for neighbors of any $b' \in B'$. It follows that $G'$ is isomorphic to the bipartite double cover of $O_{m}$, i.e. $(dO)_m$. Thus, $G'$ is regular with valency $m = \th_1$. Also $|G|=2\times \binom{2m+1}{m}=2\times \binom{2m-1}{m-1}\times \frac{(2m+1)(2m)}{(m+1)m}\geq 2|G'|$. Thus Lemma \ref{iran} applies.

\vski

The remaining infinite families are defined in terms of finite fields, so in what follows $F$ will always denote a finite field of order $q$, where $q$ is a prime power.

	\item $G$ is a Grassmann graph $J_q(n,e)$, so $G$ can be realized as the set of all $e$-dimensional subspaces of $F^n$, an $n$-dimensional vector space over $F$. Two vertices are adjacent when the intersection of their corresponding vector spaces has dimension $e-1$ (\cite[p.268]{drgraphs}). Here $\theta_1=q^2\binom{e-1}{1}_q\binom{n-e-1}{1}_q-1$ (see \cite[Thm. 9.3.3]{drgraphs}), where $\binom{m}{r}_q = \frac{(q^m-1) \ldots (q^{m-r+1}-1)}{(q^r-1) \ldots (q-1)}$ denotes the Gaussian binomial coefficient. As with the Johnson graph, we can assume $n \geq 2e$. Let us realize $F^n$ as all $n$-tuples $(a_1,\ldots, a_n)$, where $a_j \in F$. We then form a subgraph $G'$ induced on the set of all vertices corresponding to subspaces contained in the subspace $\{(a_1,\ldots, a_n) \in F^n: a_1 = 0\}$. It is easy to see that $G'$ is isomorphic to $J_q(n-1,e)$, and the valency of $G'$ is $k'=q\binom{e}{1}_q\binom{(n-1)-e}{1}_q$. Since $\binom{e}{1}_q/\binom{e-1}{1}_q\geq q$, we have $k'\geq \theta_1$. Also, $|G|=\binom{n}{e}_q=\binom{n-1}{e}_q\times \frac{q^n-1}{q^{n-e}-1}\geq 2|G'|$, so Lemma \ref{iran} applies.	

	\item $G$ is a twisted Grassmann graph $TJ_q(2e+1,e)$ (\cite{twisty}). In order to realize $G$, we consider $F^{2e+1}$ where $F$ is a filed of order $q$, which is a prime power. Fix a $(2e)$-dimensional subspace $H$ of $F^{2e+1}$. Vertices of $G$ are of two types, $(e+1)$-dimensional subspaces of $F^{2e+1}$ which are not contained in $H$ and $(e-1)$-dimensional subspaces of $H$.
This graph has the same parameters as the Grassman graph $J_q(2e+1,e)$, so in particular $\theta_1=q^2\binom{e-1}{1}_q\binom{e}{1}_q-1$. Adjacency is defined in different ways for the two types of vertices, but two vertices of the second type are adjacent if their intersection is a $(e-2)$-dimensional subspace of $F^{2e+1}$, and it follows that if we let $G'$ be induced on the set of all vertices of this type then $G'$ is isomorphic to $J_q(2e,e-1)$, which has valency $k'=q\binom{e-1}{1}_q\binom{e+1}{1}_q$. Since $\binom{e+1}{1}_q/\binom{e}{1}_q\geq q$, we have $k'\geq \theta_1$. Furthermore $|G|=\binom{2e+1}{e}_q=\frac{q^{2e+1}-1}{q^e-1} \times \binom{2e}{e-1}_q \geq 2 \times \binom{2e}{e-1}_q = 2 |G'|$, so Lemma \ref{iran} applies.

    \item $G$ is a doubled Grassman graph $dJ_q(2t+1,t)$. The method of finding a subgraph of similar type with large enough valency, which works for the other families, does not seem to work in this case. We will therefore postpone this family of graphs until Section \ref{bipbip}, where we discuss a general method which applies to bipartite graphs.

\medskip

The remaining infinite families of graphs all correspond to various sets of classical parameters $(D,b,\alpha,\beta)$. This means that their intersection arrays and spectrum are entirely determined by these four values (see \cite[Sec. 6.1 and Cor. 8.4.2]{drgraphs}). In particular, assuming $b>0$, we have

	\begin{align*}
		\theta_1 & =b^{-1}(\binom{D}{1}_b-1)(\beta-\alpha)-1 \\
		k        & =\binom{D}{1}_b\beta.
	\end{align*}

The parameters for the various families can be found on \cite[p. 24]{koolsurvey}.

	\item $G$ is a bilinear forms graph $Bil(D\times e, q)$ ($D\leq e$), with classical parameters $(D,b,\alpha,\beta)=(D,q,q-1,q^e-1)$ for $q>1$. The vertices of $G$ can be realized as the set of all bilinear maps from $F^D \times F^e$ to $F$. We may represent such a bilinear map $\phi$ by a $D \times e$ matrix $A$ with entries in $F$
such that $\phi(x,y) = x^T A y$. Two bilinear forms are adjacent if the difference of their corresponding matrices has rank 1 (see \cite[p. 280]{drgraphs}). Here $\theta_1=q^{-1}\Big(\frac{(q^{D}-1)}{q-1}-1\Big)(q^e-q) - 1 = \frac{(q^{D-1}-1)(q^e-q)}{q-1}-1$. Fix a $1$-dimensional subspace $H$ of $F^D$ and let $G'$ be induced on the set of all maps $f$ such that $f(h,v) = 0$ for all $h \in H, v \in F^{e}$. It may be checked that $G'$ is isomorphic to $Bil((D-1)\times e, q)$, which has valency $k'=\frac{(q^{D-1}-1)(q^e-1)}{q-1}> \th_1$. Furthermore it may be verified that $|G|= q^e \times |G'| \geq 2|G'|$, and so Lemma \ref{iran} applies.

	\item $G$ is an alternating forms graph $Alt(n,q)$, with classical parameters $(D,b,\alpha,\beta)=(\lfloor n/2\rfloor,q^2,q^2-1,q^m-1)$, where $m=2\lceil n/2\rceil-1$. The vertices of $G$ can be realized as the set of all bilinear maps from $F^n \times F^n$ to $F$ which are alternating (i.e. $f(x,x) = 0$ for all $x \in F^n$). If $f,g \in G$, then $f \sim g$ if $f-g$ has rank $2$ (see \cite[p. 282]{drgraphs}). Fix a $1-$dimensional subspace $H$ of $F^n$, and let $G'$ be induced on the set of all vertices corresponding to forms $f$ such that $f(x,y)=f(y,x) =0$ for all $x \in H, y \in F^n$. It may be checked that $G'$ is isomorphic to $Alt(n-1,q)$. Then
	\begin{align*}
		\theta_1 & =(\binom{\lfloor n/2\rfloor}{1}_{q^2}-1)(q^{2\lceil n/2\rceil-3}-1)-1 \\
		         & =\binom{\lfloor n/2\rfloor-1}{1}_{q^2}(q^{2\lceil n/2\rceil-1}-q^2)-1
	\end{align*}
and $k'=\binom{\lfloor (n-1)/2\rfloor}{1}_{q^2}(q^{2\lceil (n-1)/2\rceil-1}-1)$. If $n$ is even, then
	$$
		k'=\binom{\lfloor n/2\rfloor-1}{1}_{q^2}(q^{2\lceil n/2\rceil-1}-1)\geq \theta_1.
	$$
If $n$ is odd, then
	$$
		k'=\binom{\lfloor n/2\rfloor}{1}_{q^2}(q^{2\lceil n/2\rceil-3}-1)\geq \theta_1.
	$$

	Furthermore, $|G|=q^{n(n-1)/2}$ and $|G'|=q^{(n-1)(n-2)/2}$ (see \cite[p. 283]{drgraphs}), so $|G| \geq 2|G'|$, and Lemma \ref{iran} applies.

	\item $G$ is a Hermitian forms graph $Her(D,r^2)$ with classical parameters $(D,b,\alpha,\beta)=(D,-r,-r-1,-(-r)^D-1)$. Note that here that the underlying field $F$ has order $q=r^2$, with $r$ a prime power (so we require that $q$ is a prime to an even power). The vertices of $G$ can be realized as the set of all bilinear maps from $F^D \times F^D$ to $F$ which are Hermitian (i.e. $f(x,y) = f(y,x)^r$ for all $x,y \in F^D$, and note that raising $f$ to the $r$-th power is the same as applying the involutive Frobenius automorphism). If $f,g \in G$, then $f \sim g$ if $f-g$ has rank $1$ (see \cite[p. 285]{drgraphs}). Fix a $1-$dimensional subspace $H$ of $F^D$, and let $G'$ be induced on the set of all vertices corresponding to forms $f$ such that $f(x,y)=f(y,x) =0$ for all $x \in H, y \in F^D$. Then it may be checked that $G'$ is isomorphic to $Her(D-1,r^2)$ (see \cite[Prop 10.9]{grover}). Then the eigenvalues are
	\begin{align*}
		& \frac{1}{(-r)^i} (\binom{D}{1}_{-r}-\binom{i}{1}_{-r})\times(-(-r)^D-1+(r+1)\binom{i}{1}_{-r})-\binom{i}{1}_{-r} \\
		         & =\frac{1}{(-r)^i}\frac{(-r)^{2D}-(-r)^{2i}}{r+1}-\binom{i}{1}_{-r}\\
		         &=\frac{(-r)^{2D-i}-1}{r+1},
	\end{align*}
	and the second largest of these is $\theta_1 = \frac{r^{2D-2}-1}{r+1}$. Also,
	\begin{align*}
		k' & =\binom{D-1}{1}_{-r}(-(-r)^{D-1}-1)) \\
		   & =\frac{r^{2D-2}-1}{r+1}.
	\end{align*}
	so that $k' = \theta_1$. Furthermore $|G| = r^{D^2}$ and $|G'| = r^{(D-1)^2}$ (see \cite[p. 285]{drgraphs}), so $|G| \geq 2|G'|$ and Lemma \ref{iran} applies.

	\item $G$ is a quadratic forms graph $Qua(n,q)$ with classical parameters $(D,b,\alpha,\beta)=(\lfloor \frac{n+1}{2}\rfloor, q^2,q^2-1,q^m-1)$, where $m=2\lfloor \frac{n}{2}\rfloor +1$. A quadratic form $\ga$ on $F^n$ is a map from $F^n$ to $F$ satisfying $\ga(\la x) = \la^2 \ga(x)$ for any $\la \in F, x \in F^n$ and such that $f(x,y) = \ga(x+y)-\ga(x)-\ga(y)$ is a symmetric bilinear form. The vertices of the graph $G$ can be realized as the set of all quadratic forms from $F^n$ to $F$. If $\ga,\dd \in G$, then $\ga \sim \dd$ if $\ga-\dd$ has rank $1$ or $2$ (see \cite[p. 290]{drgraphs}). Fix a $1-$dimensional subspace $H$ of $F^D$, and let $G'$ be induced on the set of all vertices corresponding to quadratic forms that vanish on $H$. Then it may be checked that $G'$ is isomorphic to $Qua(n-1,q)$. Furthermore, $Qua(n,q)$ and $Qua(n-1,q)$ have the same parameter sets as $Alt(n+1,q)$ and $Alt(n,q)$, respectively, so the calculations done earlier for the alternating forms graph apply here and show that $k' \geq \theta_1$ and $|G| \geq 2|G'|$. Lemma \ref{iran} therefore applies.

	\item $G$ is a dual polar graph, with classical parameters $(D,b,\alpha,\beta)=(D,q,0,q^e)$ (for some $e\in\{0,\frac{1}{2},1,\frac{3}{2},2\}$). In order to describe $G$, let $V$ be a vector space over $F$ with an associated type of form as specified on \cite[p. 274]{drgraphs}; note that the dimension of $V$ is $2D, 2D+1,$ or $2D+2$. The vertices are realized as the maximal isotropic subspaces of $V$({\it isotropic} means that a form of the given type vanishes there), which are necessarily of dimension $D$, and two vertices are adjacent if the intersection of the corresponding subspaces has dimension $D-1$. Let $H$ be an isotropic 1-dimensional space and $W$ the orthogonal complement of $H$. Now let $P$ be a hyperplane in $W$, and let the graph $G'$ be induced on the vertices corresponding to maximal isotropic spaces in $P$. It may be checked then that $G'$ is isomorphic to the same type of dual polar graph with parameters $(D-1,q,0,q^e)$. Then we have $\theta_1=\frac{1}{q}(\binom{D}{1}_q-1)q^e-1= \binom{D-1}{1}_q q^e-1$ and $k'=\binom{D-1}{1}_q q^e$ so $k'\geq \theta_1$. Also $|G|=\prod_{i=0}^{D-1}(q^{i+e}+1)$ and $|G'|=\prod_{i=0}^{D-2}(q^{i+e}+1)$ by \cite[Lem. 9.4.1]{drgraphs}, so that $|G| \geq 2|G'|$. Thus Lemma \ref{iran} applies.

\vski

A Hemmeter graph $G$ is the extended bipartite double of the dual polar graph on $\mathcal C_{D}(q)$, and has the same parameter set as the dual polar graph on $\mathcal D_{D+1}(q)$ (see \cite[p.279]{drgraphs}). Since we saw above that the dual polar graph on $\mathcal C_{D}(q)$ contains the dual polar graph on $\mathcal C_{D-1}(q)$ as an induced subgraph, the extended bipartite double of this induced subgraph will be a subgraph of $G$ with the same parameters as the dual polar graph on $\mathcal D_{D}(q)$. Therefore the calculations above for the dual polar graphs also prove the conjecture for the Hemmeter graphs.
	
\item $G$ is a half dual polar graph $D_{n,n}(q)$. $G$ is the halved graph of the dual polar graph on $\mathcal D_{D}(q)$ (see \cite[p.278]{drgraphs}). Here we have classical parameters $(D,b,\alpha,\beta)=(\lfloor \frac{n}{2}\rfloor,q^2,q^2+q,\frac{q^{m+1}-1}{q-1}-1)$, where $m=2\lceil \frac{n}{2}\rceil-1$. We may take as $G'$ the halved graph of the subgraph defined for the dual polar graph, and the relation $|G| \geq 2|G'|$ is preserved by halving. Furthermore we have $\theta_1=q^3\binom{D-2}{2}_{q^2}-1$ and $k'=q\binom{D-1}{2}_{q^2}$ (see \cite[Thm. 9.4.8]{drgraphs}). It is straightforward check that $q\binom{D-1}{2}_{q^2} > q^3\binom{D-2}{2}_{q^2}$, and thereby to verify $k' > \theta_1$. Lemma \ref{iran} therefore applies.

\vski

The Ustimenko graphs are the halved graphs of the Hemmeter graphs, and if $G$ is such a graph we may take as $G'$ the halved graph of the the subgraph defined above for the Hemmeter graphs; again the relationship $|G| \geq 2|G'|$ persists, and the shared parameters with the dual polar graphs on $\mathcal D_{D}(q)$ prove the conjecture in this case.
\end{enumerate}

{\bf Remark:} The subgraphs $G'$ found for the infinite families above are essentially the {\it descendents} discussed in the work \cite{tanakadescend}, which in turn builds upon a result from \cite{kooldescend}.

\section{Strongly-regular graphs} \label{srg}

A {\it strongly regular graph} is simply a distance-regular graph of diameter 2. The parameters for such graphs are commonly given with the notation $srg(v, k, \la, \mu)$, which means that $G$ is a strongly regular graph with $v$ vertices, valency $k$, and $\la, \mu$ are the same quantities that we have been calling $a_1,c_2$, respectively. To avoid confusion with the Laplacian eigenvalues $\la_i$, and to keep the notation consistent throughout the paper, we will use the $a_i, b_i, c_i$ notation for strongly regular graphs as well (this also makes sense because of the importance of $b_1=k-\la-1$ in our analysis below). Thus, we may present the parameters for a strongly regular graph $G$ by stating that $G$ is $srg(n,k,a_1,c_2)$; note that $k,a_1,c_2$ uniquely determine all other values in the intersection array. In this and subsequent sections we will describe a graph as "OK" whenever we can exhibit a set which gives an upper bound for $h_G$ which is at most $\la_1$.

\begin{proposition} \label{cherry}
Every strongly regular graph is OK.
\end{proposition}

{\bf Proof:} We will prove this through a short collection of lemmas. We note first that a logical disconnecting set of edges to consider in a strongly regular graph is the set of edges connecting $\Ga_1(x)$ and $\Ga_2(x)$ for a given vertex $x$. This set of $kb_1 = |\Ga_2(x)|c_2$ edges separates the graph into the disjoint sets $\{x\} \cup \Ga_1(x)$ and $\Ga_2(x)$. If $|\{x\} \cup \Ga_1(x)|\leq|\Ga_2(x)|$, then we obtain $\frac{kb_1}{k(k+1)} = \frac{b_1}{k+1}$ as an upper bound for $h_\GGG$, whereas if $|\{x\} \cup \Ga_1(x)|\geq|\Ga_2(x)|$ then we obtain an upper bound of $\frac{|\Ga_2(x)|c_2}{|\Ga_2(x)|k} = \frac{c_2}{k}$. This shows easily that $h_\GGG \leq \max\(\frac{b_1}{k+1}, \frac{c_2}{k}\)$, and we obtain

\begin{lemma} \label{kristi}
If $G$ is a strongly regular graph and $\max\(\frac{b_1}{k+1}, \frac{c_2}{k}\) \leq \la_1$, then $G$ is OK.
\end{lemma}

This turns out to be sufficient in many cases, though not in all.  Note that every strongly regular graph has precisely 3 distinct eigenvalues: $k > \th_1>\th_2$. These eigenvalues must be integers unless the graph in question is a conference graph. It therefore simplifies matters to dispose of the conference graphs first.

\begin{lemma} \label{}
Conference graphs are OK.
\end{lemma}

{\bf Proof:} These are graphs of the form $srg(v,\frac{v-1}{2},\frac{v-5}{4},\frac{v-1}{4})$. Here it may be checked that $\max\(\frac{b_1}{k+1}, \frac{c_2}{k}\) = \frac{c_2}{k} = \frac{1}{2}$. The eigenvalues of the graph are $\frac{v-1}{2}, \frac{\sqrt{v}-1}{2},$ and $\frac{-\sqrt{v}-1}{2}$, which gives $\la_1 = \frac{v-\sqrt{v}}{v-1}$, and this quantity is greater than $\frac{1}{2}$ for $v > 1$ and $v \equiv 1 (\mbox{mod }4)$. \qed

Henceforth we assume that the eigenvalues of the graphs in question are integers. In order to deal with the remaining cases we need the following result.

\begin{lemma} \label{emily}
Let $G$ be a strongly regular graph on $v$ vertices. If $v=2t$ or $v=2t+1$, and $tk-(2t+1)\th_1 +(t+1)\th_2 \geq 0$, then $G$ is OK.




\end{lemma}

In order to prove this, we need another lemma.

\begin{lemma} \label{venus}
Let $G$ be a regular graph with valency $k$, $v$ vertices and smallest eigenvalue $\theta_{\min}$.
Then the following holds:
	\begin{itemize}
		\item If $v=2t$, then $h_G\leq \frac{(k-\theta_{\min})}{2k} \leq \frac{(t+1)(k-\theta_{\min})}{(2t+1)k}$;
		\item If $v=2t+1$, then $h_G\leq \frac{(t+1)(k-\theta_{\min})}{(2t+1)k}$.
	\end{itemize}
\end{lemma}

{\bf Proof:} If $v=2t$ is even, then we take any partition $\pi=\{S,S^c\}$ of $V(G)$ with $|S|=|S^c|=t$. If we set $\aaa=\frac{E[S,S^c]}{t}$ to be the average number of vertices in $S^c$ that each point in $S$ is adjacent to, then following for instance the methods found in \cite[p. 596]{eigint} we have the following quotient matrix:
	$$
	\begin{pmatrix}
		k-\alpha & \alpha\\
		\alpha & k-\alpha
	\end{pmatrix}.
	$$
	This matrix has eigenvalues $k$ and $k-2\aaa$, and by the interlacing described in \cite[Cor. 2.3]{eigint} we see that $k-2\aaa \geq \theta_{\min}$. It follows that $h_G\leq \frac{\alpha}{k}\leq \frac{k-\theta_{\min}}{2k}$.

\vski

	The method for the case $v=2t+1$ is essentially identical. We again take any partition $\pi=\{S,S^c\}$ with $|S|=t$, and arguing as before gives the quotient matrix
	$$
	\begin{pmatrix}
		k-\alpha & \alpha\\
		\beta & k-\beta
	\end{pmatrix},
	$$ where the relation $t\alpha=(t+1)\beta$ holds. This matrix has eigenvalues $k$ and $k-\aaa-\beta$, and the interlacing argument again gives $k-\aaa-\beta \geq \theta_{\min}$. It follows that $h_G\leq \frac{\alpha}{k}\leq \frac{(t+1)(k-\theta_{\min})}{(2t+1)k}$. \qed

Lemma \ref{emily} follows easily from Lemma \ref{venus}, as the conditions of Lemma \ref{emily} can be rearranged to $\frac{(t+1)(k-\th_2)}{(2t+1)k} \leq \frac{k-\th_1}{k} = \la_1$, and the result follows. In order to complete the proof of Proposition \ref{cherry}, we will show that either Lemma \ref{kristi} or Lemma \ref{emily} can be applied in every case.

\vski

In order to accomplish this, let us assume that the conditions of Lemma \ref{kristi} do not hold; that is, that $\max\(\frac{b_1}{k+1}, \frac{c_2}{k}\) > \la_1$. It is a standard fact that $-\th_1\th_2 = k-c_2$ (see for instance \cite[Thm. 1.3.1 (iii)]{drgraphs}); inserting $\th_1 = \frac{k-c_2}{-\th_2}$ into $\la_1 = \frac{k-\th_1}{k}$ yields $\la_1 = \frac{(-\th_2-1)k+c_2}{-\th_2k}$. This is at least $\frac{c_2}{k}$ whenever $(-\th_2-1)k \geq (-\th_2-1)c_2$, and this is always true since $k \geq c_2$. Thus, we may assume $\frac{b_1}{k+1} > \la_1 =\frac{k-\th_1}{k}$, and this implies $b_1+\th_1 \geq k+1$ since $b_1$ and $\th_1$ are integers. 

\vski

Let $t = \lfloor \frac{v}{2} \rfloor$. By Lemma \ref{emily} we are now required to prove $tk \geq (2t+1)\th_1 -(t+1)\th_2$, and as before we may assume that $b_1+\th_1 \geq k+1$. Suppose that $\th_1 \geq 3, \th_2 \leq -3$. Then $|\th_1|,|\th_2| \leq \frac{k-c_2}{3}$, and it will suffice now to prove that $tk \geq (2t+1)\frac{k-c_2}{3} + (t+1)\frac{k-c_2}{3} = tk -tc_2 + \frac{2(k-c_2)}{3}$. This will hold whenever $\frac{3}{2}tc_2 \geq k$: if $c_2 \geq 2$ then we are reduced to showing $3t \geq k$, which is clearly true since $3t \geq v$; while if $c_2 = 1$ then we must have $k \leq t$, since then $2t+1 \geq v = 1 + |\Ga_1| + |\Ga_2| = 1 + |\Ga_1| + |\Ga_1|b_1 \geq 2k+1$, so that $\frac{3}{2}tc_2 \geq k$ holds as well. Thus, at least one of $\th_1, \th_2$ must have modulus at most 2.

\begin{itemize} \label{}

\item If $\th_1 =2$, then $\th_2 = \frac{-(k-c_2)}{2}$, and we must show $tk \geq 4t+2 +(t+1)\frac{k-c_2}{2}$, or equivalently $(t-1)k + (t+1)c_2 \geq 8t + 4$. Here $b_1+\th_1 \geq k+1$ implies $b_1 = k-1$ and $a_1 = 0$. Using $a_1-c_2 = \th_1 + \th_2$ (\cite[Thm. 1.3.1(iii)]{drgraphs}) we see that $-c_2 = \th_1 + \th_2 = 2 - \frac{k-c_2}{2}$, which implies $k = 4+3c_2$. If $c_2 \geq 2$, so that $k \geq 10$, then it is easy to see that $(t-1)k + (t+1)c_2 \geq 8t + 4$ holds since $t \geq 3$. If $c_2 = 1$, then the graph in question, with $k=7, b_1=6, c_2 =1$, has 50 vertices. In this case, $\th_1 = 2$ and $\th_2=-3$, and we can apply the bound $h_G\leq \frac{(k-\theta_{\min})}{2k}$ in Lemma \ref{venus} to obtain $h_G \leq \frac{5}{7} = \la_1$.

\item If $\th_1 = 1$, then $b_1+\th_1 \geq k+1$ implies $b_1 \geq k$, a contradiction.

\item If $\th_1 =0$ then the proposition holds trivially.

\item $\th_2 = -1$ occurs only for the complete graph.

\item Suppose $\th_2 = -2$. In this case $\th_1 = \frac{k-c_2}{2}$. 
According to \cite[Thm. 3.12.4]{drgraphs}, the only strongly regular graphs with smallest eigenvalue $-2$ are the triangular graph $T(m) (m \geq 5)$, the lattice graph $L_2(m) (m \geq 3)$, the complete multipartite graph $K_{m \times 2}$, the Petersen graph, the Clebsch graph, the Schl\"afli graph, the Shrikhande graph, and the Chang graphs. In this list we must only address the graphs with $\th_1 \geq 3$, because the others have been dealt with above. This leaves us only $L_2(m) (m \geq 5)$, $T(m) (m \geq 7)$, the Schl\"afli graph, and the Chang graphs.
Note that in these cases $\la_1 = \frac{k-\th_1}{k} = \frac{k+c_2}{2k}$. The Schl\"afli graph is $srg(27,16,10,8)$, so $\max\(\frac{b_1}{k+1}, \frac{c_2}{k}\) = \frac{c_2}{k} = \frac{1}{2} < \la_1$. The lattice graph $L_2(m)$ is $srg(m^2,2m-2,m-2,2)$ with $m \geq 3$, and $\max\(\frac{b_1}{k+1}, \frac{c_2}{k}\) = \frac{b_1}{k+1} =\frac{m-1}{2m-1} < \frac{1}{2} < \la_1$. The triangular graph $T(m)$ is $srg(m(m-1)/2,2m-4,m-2,4)$ with $m \geq 5$, and $\max\(\frac{b_1}{k+1}, \frac{c_2}{k}\) = \frac{b_1}{k+1}$ for $m \geq 8$, and $\frac{b_1}{k+1} = \frac{m-3}{2m-3} < \frac{1}{2} < \la_1$.; for $m = 5,6,7$ we have $\max\(\frac{b_1}{k+1}, \frac{c_2}{k}\) = \frac{c_2}{k} = \frac{4}{k} < \frac{k+4}{2k} = \la_1$. Finally, the Chang graphs all have the same parameters as $T(8)$ and are therefore covered by the same argument.

\end{itemize}

This completes the proof of Proposition \ref{cherry}.

\section{Several families of distance-regular graphs with diameter 3} \label{diam3}

In this section we prove the conjecture for a number of families of distance-regular graphs with diameter 3. 

\subsection{Bipartite graphs} \label{bipbip}

Here we study a construction that allows us to prove that a number of bipartite distance-regular graphs are OK. First, two necessary lemmas.

\begin{lemma} \label{mar}
Let $G$ be any graph. Suppose $A$ and $B$ are disjoint subsets of $V(G)$, and $|B| = r$. Then, for any positive integer $r' \leq r$ there is a subset $B'$ of $B$ with $|B'| = r'$ such that $E[A,B'] \geq \frac{r'}{r}E[A,B]$.
\end{lemma}

{\bf Proof:} We have

\begin{equation} \label{}
\sum_{\substack{B' \subseteq B \\ |B'|=r'}} E[A,B'] = \sum_{b \in B} {r-1 \choose r'-1} E[A,\{b\}] = {r-1 \choose r'-1} \sum_{b \in B} E[A,\{b\}] = {r-1 \choose r'-1} E[A,B].
\end{equation}

Thus the average size of $E[A,B']$ over all such $B'$ is $\frac{{r-1 \choose r'-1}}{{r \choose r'} } E[A,B] = \frac{r'}{r}E[A,B]$, and the result follows. \qed

\begin{lemma} \label{tie}
Let $G$ be a connected bipartite regular graph with $2r$ vertices and valency $k$.
	\begin{itemize}
		\item If $r=2m$ is even, then $h_G\leq \frac{1}{2}$;
		\item If $r=2m+1$ is odd, then $h_G\leq \frac{1}{2} + \frac{1}{8(m+\frac{1}{2})^2} = \frac{1}{2} + \frac{1}{2r^2}$.
	\end{itemize}
\end{lemma}

{\bf Proof:} Let the two bipartite components be denoted $A$ and $B$. Suppose first that $r=|A|=|B|$ is even. If we take a subset $A_1\subset A$ with $|A_1|=m=\frac{r}{2}$, then by Lemma \ref{mar} there exists a subset $B_1\subset B$ with $|B_1|=m$ such that $E[A_1,B_1]\geq \frac{|B_1|}{|B|} E[A_1,B]=\frac{km}{2}$. If we set $S=A_1\cup B_1$, we see that the subgraph $G'$ induced on $S$ has average valency $k' = \frac{2E[A_1,B_1]}{|A_1|+|B_1|}\geq \frac{k}{2}$ and satisfies $|G'|=|G|/2$, and the result follows since $h_G \leq \frac{k-k'}{k}$ (see the proof of Lemma \ref{iran}). Now suppose that $r=|A|=|B|$ is odd. We may take a subset $A_1\subset A$ with $|A_1|=m=\frac{r-1}{2}$, and then by Lemma \ref{mar} there exists a subset $B_1\subset B$ with $|B_1|=m+1$ such that $E[A_1,B_1]\geq \frac{|B_1|}{|B|} E[A_1,B]=\frac{km(m+1)}{2m+1}$. If we set $S=A_1\cup B_1$, we see that the subgraph $G'$ induced on $S$ has average valency $k' = \frac{2E[A_1,B_1]}{|A_1|+|B_1|}\geq \frac{2 km(m+1)}{(2m+1)^2} = \frac{k((m+1/2)^2-1/4)}{2(m+1/2)^2}$ and satisfies $|G'|=|G|/2$. Then $h_G \leq \frac{k-k'}{k} \leq 1-\frac{((m+1/2)^2-1/4)}{2(m+1/2)^2}$ and the result follows. \qed

Armed with Lemma \ref{tie}, let us begin by discussing with the doubled Grassmann graphs (Family 11 in Section \ref{infam}). If $G=dJ_q(2t+1, t)$ is such a graph, then $G$ has as halved graphs the Grassmann graphs $J_q(2t+1, t)$ (see \cite[Th. 9.3.11]{drgraphs}). It follows that the valency $k$ and second largest eigenvalue $\theta_1$ of $dJ_q(2t+1, t)$ satisfy $k = { t+1 \choose 1}_q = \frac{q^{t+1}-1}{q-1}$ and $(\theta_1^2 - k)/c_2 = \tilde{\theta}_1$ where $\tilde{\theta}_1$ is the second largest eigenvalue of $J_q(2t+1, t)$. As $c_2 =1$ (for $dJ_q(2t+1, t)$) we obtain after a calculation $\th_1^2= q^2 { t-1 \choose 1}_q { t \choose 1}_q - 1 + { t+1 \choose 1}_q= q { t \choose 1}^2_q$, hence $\th_1 = \sqrt{q} { t \choose 1}_q$. Thus,

\be \lll{swim}
\la_1 = \frac{k-\th_1}{k} = \frac{q^{t+1}-1 - \sqrt{q} (q^{t}-1)}{q^{t+1}-1} \geq 1 - q^{-(t+1)} - q^{-1/2} + q^{-(t+1/2)} \geq 1- \frac{1}{\sqrt{q}}.
\ee

For $q \geq 5$ we obtain $\la_1 \geq .55$, while even assuming $r$ is odd and $m\geq 2$ Lemma \ref{tie} gives an upper bound of $.52$ for $h_G$, so that $G$ is OK for $q \geq 5$. If $q=4$ and $r$ is even, then again Lemma \ref{tie} applies, although in fact $r$ will always be odd for $q=4$. This is because the halved graph $J_q(2t+1, t)$ has

\be \lll{breaststroke}
r = { 2t+1 \choose t}_4 = \frac{(4^{2t+1} - 1)(4^{2t} - 1)\ldots (4^{t+2} - 1)}{(4^{1} - 1) (4^{2} - 1) \ldots (4^{t} - 1)}
\ee

vertices, and this must always be odd. Nevertheless, Lemma \ref{tie} does suffice here, since \rrr{swim} does show that $\la_1 \geq \frac{1}{2} + (\frac{1}{4^{t+1/2}} - \frac{1}{4^{t+1}}) = \frac{1}{2} + \frac{1}{4^{t+1}}$, and it is straightforward to verify that $\frac{1}{2r^2} < \frac{1}{4^{t+1}}$ using \rrr{breaststroke}. The lemma does not work for $q=2,3$ (except when $t=1$), since in those cases for $t>1$ we have $\la_1 < \frac{1}{2}$. These cases remain open. We summarize with a proposition.

\begin{proposition} \label{}
For $q \geq 4$, the doubled Grassmann graph $dJ_q(2t+1, t)$ is OK.
\end{proposition}

\vski

By the well-known Feit-Higman Theorem \cite[Thm 6.5.1]{drgraphs}, a bipartite distance-regular graph with diameter $D$ and valency $k \geq 3$
with $c_{D-1}=1$ only exists if $D \in \{2,3,4,6\}$. Such graphs with $D=2$ are strongly regular, and are therefore OK by Proposition \ref{cherry}, while those with diameter $D=3$ are covered by Proposition \ref{bip3}. The case $D=4$ correspond to the incidence graphs of generalized quadrangles of order $q$. These graphs have intersection array $\{q+1, q, q, q; 1, 1, 1, q+1\}$ and $\theta_1 = \sqrt{2q}$, and they exist whenever $q$ is a prime power. The cases $q=2,3$ will be shown to be OK in Section \ref{smallval}. For larger $q$, we have the following proposition.

\begin{proposition} \label{}
For $q \geq 7$, if $G$ is the incidence graph of a $GQ(q, q)$, then $G$ is OK.
\end{proposition}

{\bf Proof:} $G$ has $v=2(q^2+1)(q+1)$ vertices, and if $q$ is odd then $v$ is a multiple of 4, so Lemma \ref{tie} implies $h_G \leq \frac{1}{2}$. Since $\la_1 = \frac{q+1-\sqrt{2q}}{q+1}$ is increasing in $q$ and equal to $\frac{8-\sqrt{14}}{8}>\frac{1}{2}$ for $q=7$, we see that $G$ is OK for odd $q \geq 7$. If $q \geq 8$ is a power of 2, then $\la_1 \geq \frac{9-\sqrt{16}}{9} = \frac{5}{9}$, while by Lemma \ref{tie} we have $h_G \leq \frac{1}{2} + \frac{2}{v^2} \leq \frac{1}{2} + \frac{2}{1170^2}< \frac{5}{9}$, so $G$ is OK. \qed

{\bf Remark:} The cases $q=4,5$ remain open.

\vski

We may also address the incidence graphs of generalized hexagons of order $(q, q)$. These have intersection array $\{q+1, q, q, q, q, q; 1, 1, 1, 1, 1, q+1\}$ and $\theta_1 = \sqrt{3q}$, and again exist whenever $q$ is a prime power. The case $q=2$ is Tutte's 12-cage, and we will show this to be OK in Section \ref{smallval}. As indicated also in that section, we have not been able to prove that the graph is OK for $q=3$. For larger $q$ we have the following.

\begin{proposition} \label{}
For $q \geq 11 $ the incidence graph of a $GH( q, q)$ is OK.
\end{proposition}

{\bf Proof:} $G$ has $v=2(q^4+q^2+1)(q+1)$ vertices, and if $q$ is odd then $v$ is a multiple of 4, so Lemma \ref{tie} implies $h_G \leq \frac{1}{2}$. Since $\la_1 = \frac{q+1-\sqrt{3q}}{q+1}$ is increasing in $q$ and equal to $\frac{12-\sqrt{33}}{12}>\frac{1}{2}$ for $q=11$, we see that $G$ is OK for odd $q \geq 11$. If $q \geq 16$ is a power of 2, then $\la_1 \geq \frac{17-\sqrt{48}}{17} > \frac{10}{17}$, while by Lemma \ref{tie} we have $h_G \leq \frac{1}{2} + \frac{2}{v^2} \leq \frac{1}{2} + \frac{2}{2236962^2}< \frac{10}{17}$, so $G$ is OK. \qed

{\bf Remark:} The cases $q=3,4,5,7,8,$ and $9$ remain open.

\vski


We also have the following.

\begin{proposition} \label{bip3}
Bipartite distance-regular graphs of diameter 3 are OK.
\end{proposition}

{\bf Proof:} By \cite[p. 432]{drgraphs}, such graphs satisfy $\th_1 = \sqrt{k-c_2} \leq \sqrt{k-1}$, and thus $\la_1 \geq \frac{k-\sqrt{k-1}}{k}$. For $k \geq 4$, this gives a bound of $\la_1 \geq \frac{28}{50}$. Note that for diameter 3 we have $|G| \geq 2k+2$, since $|\Ga_0| + |\Ga_2| = |\Ga_1|+|\Ga_3| \geq k+1$. Thus, $k \geq 4$ and $|G| = 4m$ or $4m+2$ implies $m \geq 2$, so that Lemma \ref{tie} implies $h_G \leq \frac{26}{50}$, and the result follows. Note that the case $k=3$ is covered in Section \ref{smallval}. \qed

\subsection{Antipodal distance-regular graphs with diameter three}

\begin{proposition} \label{cosplay}
	Antipodal distance-regular graphs with diameter 3 are OK.
\end{proposition}

{\bf Remark:} Note that this class of graphs includes the Taylor graphs.

\vski

{\bf Proof:} Let $G$ be such a graph. Let $t=\lfloor \frac{k+1}{2} \rfloor$, and let us first suppose that $\th_1 \leq \frac{t}{k}b_1$. Fix $x_0 \in V(G)$, and write $\Ga_3(x_0) = \{x_1, \ldots, x_r\}$. Note that $d(x_i,x_j) = 3$ for any distinct $i,j \in \{0, \ldots, r\}$. Let $X_j = \Ga_1(x_j)$ for $j \in \{0, \ldots, r\}$.  Choose an arbitrary set $A_0\subset X_0$ with $|A_0|=t$. Then, by Lemma \ref{mar}, there exists a set $A_1\subset X_1$ with $|A_1|=t$, such that $E[A_0,A_1]\geq \frac{|A_1|}{|X_1|}E[A_0,X_1]= (\frac{t}{k})tc_2$. If we set $B_1=A_0\cup A_1$, then $\frac{E[B_1,B_1]}{|B_1|} \geq \frac{2E[A_0,A_1]}{|B_1|}\geq \frac{2(\frac{t}{k})tc_2}{2t} = \frac{tc_2}{k}$. Assume now that for $j\leq r-1$ we have defined $B_j=\cup_{i=0}^j A_i$ with $A_i\subset X_i$ and $|A_i|=t$ for all $i$, and $\frac{E[B_j,B_j]}{|B_j|}\geq \frac{jt}{k}c_2$. Then there exists a subset $A_{j+1}\subset X_{j+1}$ with $|A_{j+1}|=t$, such that $E[B_j,A_{j+1}]\geq \frac{|A_{j+1}|}{|X_{j+1}|}E[B_j,X_{j+1}]\geq (\frac{t}{k}) (j+1)tc_2$ (using here $|B_j| = (j+1)t$). If we take $B_{j+1}=B_j\cup A_{j+1}$, then $\frac{E[B_{j+1},B_{j+1}]}{|B_{j+1}|}\geq\frac{E[B_j,B_j]+2E[B_j,A_{j+1}]}{|B_{j+1}|}\geq \frac{(j+1)t(\frac{jt}{k}c_2)+2(\frac{t}{k})(j+1)tc_2}{(j+2)t}=\frac{(j+1)t}{k}c_2$. Continuing in this manner, we may inductively define a set $B_r=\cup_{i=0}^r A_i$ with $A_i\subset X_i$ and $|A_i|=t$, such that $\frac{E[B_r,B_r]}{|B_r|}\geq \frac{rt}{k}c_2$. Since $G$ is antipodal, we have $b_2 = 1$ and $c_3 = k$, thus $r = |\Ga_3| = \frac{b_0 b_1 b_2}{c_1 c_2 c_3} = \frac{b_1}{c_2}$, and $rc_2 = b_1$. Therefore our set $B_r$ satisfies $\frac{E[B_r,B_r]}{|B_r|}\geq\frac{t}{k}b_1$; in other words if we let $G'$ be induced on $B_r$, then the average valency of $G'$ is at least $\frac{t}{k}b_1 \geq \th_1$. Furthermore $|G'| = (r+1)t$, but $|G| = \sum_{j=0}^{r}(|X_j|+1) \geq 2(r+1)t$, and thus Lemma \ref{iran} applies.

\vski

Now let us suppose that $\th_1 \leq a_1 + 1$. If we let $G'$ be induced on $\{x\} \cup \Ga_1(x)$ then clearly $2|G'| \leq |G|$ and $G'$ has average valency $\frac{E[G',G']}{|G'|} = \frac{k(a_1+1)+k}{k+1} = \frac{k}{k+1}(a_1+2) > a_1+1$, using $k+1 > a_1+2$. Thus, Lemma \ref{iran} applies, and $G$ is OK.

\vski

We may therefore assume $\th_1 > \max(\frac{t}{k}b_1, a_1+1)$. In this case $k = b_1 + a_1 + 1 < 2 \th_1 + \th_1$, and thus $\th_1 > \frac{k}{3}$. The bipartite graphs have been covered in Proposition \ref{bip3}, and for the remaining ones by \cite[p. 431]{drgraphs} we have

\begin{equation*} \label{}
\begin{gathered}
		\theta_1+\theta_3=a_1-c_2,\\
		\theta_1\theta_3=-k,
\end{gathered}
\end{equation*}

with either $\theta_1$ and $\theta_3$ both integers, or $\theta_1=-\theta_3=\sqrt{k}$ and $a_1=c_2$. If $\theta_1$ and $\theta_3$ are integers then since $\th_1>a_1+1$ we have $a_1-c_2=\theta_1+\theta_3>a_1+1+\theta_3$, which implies $\theta_3<-c_2-1\leq -2$. On the other hand, $\theta_1>\frac{k}{3}$ and $\theta_1\theta_3=-k$ implies $\theta_3\geq -2$, a contradiction.

\vski

In the case $\theta_1=-\theta_3=\sqrt{k}$, we see that $\sqrt{k}=\theta_1>\frac{k}{3}$ implies $k\leq 8$, and then $3>\theta_1>a_1+1$ implies $a_1=c_2=1$. If $k=8$, we have $\th_1 >\frac{t}{k}b_1 = \frac{b_1}{2}=3>\sqrt{k}$, a contradiction. If $k=7$, we have $\th_1 \geq \frac{t}{k}b_1=\frac{k+1}{2k}b_1=\frac{20}{7}>\sqrt{k}$, again a contradiction. If $k \leq 6$, then as the subgraph $G'$ induced on $\{x\} \cup \Ga_1(x)$ has average valency $\frac{k(a_1+1)+k}{k+1} = \frac{3k}{k+1} > \sqrt{k} = \th_1$, $G$ is OK. This covers all cases, and the proof is complete. \qed

\subsection{Shilla distance-regular graphs}

A distance-regular graph of diameter 3 is called {\it Shilla} if it satisfies $\th_1 = a_3$ (see \cite{shilla}). In this case, $\la_1 = \frac{k-\th_1}{k} = \frac{c_3}{k}$, while if we fix $x\in V(G)$ and let $S=\Ga_3(x)$, then $h_G \leq \frac{E[S,S^c]}{k|S|} = \frac{c_3|\Ga_3(x)|}{k|\Ga_3(x)|} = \frac{c_3}{k}=\la_1$, and we obtain the following.

\begin{proposition} \label{}
If $G$ is a Shilla distance-regular graph, then $G$ is OK provided that $|\Ga_3| \leq \frac{|G|}{2}$.
\end{proposition}

We remark that clearly $|\Ga_3| \leq \frac{|G|}{2}$ whenever $c_3 \geq b_2$, since in that case $|\Ga_3| \leq |\Ga_2|$, and in fact if we calculate $|G|$ precisely we find that $b_1b_2 \leq c_3(c_2+b_1)$ is sufficient. Most (though not all) of known Shilla distance-regular graphs satisfy this (see \cite{shilla} for a number of examples); an exception is the Odd graph $O_4$, which we proved to be OK by other methods in Section \ref{infam} (see also Section \ref{smallval}). Proving that all Shilla distance-regular graphs are OK remains an open problem.



\section{Graphs with small valency} \label{smallval}

In this section we dispose of most of the distance-regular graphs with valency 3 and 4, utilizing the fact that all such graphs are known. We begin with a simple observation.

\begin{lemma} \label{}
Suppose $G$ is a regular graph on $n$ vertices with valency $k \geq 3$, diameter $D \geq 3$, and girth $g$. Then $g \leq \frac{n}{2}$.
\end{lemma}

{\bf Proof:} Let $S$ be a cycle of length $g$. Each vertex in $S$ must be adjacent to precisely two other vertices in $S$ (if more than two, then we could form a shorter cycle), so since $k \geq 3$ each vertex is adjacent to at least one vertex in $S^c$. So, for every $z \in S$ we may let $\phi(z)$ be an adjacent vertex in $S^c$. If $\phi$ is injective then $|S^c| \geq |S|$ and the result follows, so we may assume that $\phi$ is not injective, and therefore there are $x, y \in S$ such that $\phi(x)=\phi(y) =: z$. The shortest path $P$ in $S$ from $x$ to $y$ can have length at most 2, since if it were longer then we could shorten the cycle by replacing $P$ by the path $\{x,z,y\}$. But, in that case, $P \cup \{z\}$ forms a cycle of length at most 4, and thus $g \leq 4$. However, the conditions on $G$ imply $n \geq 8$: given any vertex $x \in V(G)$ we can find $y \in V(G)$ with $d(x,y)\geq 3$, and then $n \geq |\{x\} \cup \Gamma_1(x) \cup \Gamma_1(y) \cup \{y\}| \geq 8$. The result follows. \qed

This lemma shows that we may always choose a cycle as $S$ in order to obtain an upper bound for $h_G$: each vertex in $S$ will have at most $k-2$ neighbors outside $S$, so that $E[S,S^c] \leq (k-2)|S|$, and we obtain $h_G \leq \frac{k-2}{k}$. This is clearly not a particularly strong bound for large valency, but it suffices to prove the conjecture for many graphs of small valency. We therefore isolate it as a lemma.

\begin{lemma} \label{arms}
Let $G$ be a distance-regular graph with valency $k \geq 3$ and diameter $D \geq 3$. Then $h_G \leq \frac{k-2}{k}$.
\end{lemma}

We begin with an analysis of the valency three graphs. According to \cite[Thm. 7.5.1]{drgraphs}, the only graphs with valency $k=3$ and diameter $D \geq 3$ are the following:

\begin{itemize} \label{}

\item The cube, with intersection array $\{3,2,1;1,2,3\}$. $\th_1 = 1$ and thus $\la_1 = \frac{2}{3}$, and Lemma \ref{arms} applies (in fact, Theorem \ref{mass} implies $h_G = \frac{1}{3}$).

\item The Heawood graph, with intersection array $\{3,2,2;1,1,3\}$. $\th_1 = \sqrt{2}$ and thus $\la_1 = \frac{3-\sqrt{2}}{3} \approx .53$, and Lemma \ref{arms} applies.

\item The Pappus graph, with intersection array $\{3,2,2,1;1,1,2,3\}$. $\th_1 = \sqrt{3}$ and thus $\la_1 = \frac{3-\sqrt{3}}{3} \approx .42$, and Lemma \ref{arms} applies.

\item The Coxeter graph, with intersection array $\{3,2,2,1;1,1,1,2\}$. $\th_1 = 2$ so that $\la_1 = \frac{1}{3}$, and Lemma \ref{arms} applies.

\item Tutte's 8-cage, with intersection array $\{3,2,2,2;1,1,1,3\}$. $\th_1 = 2$ so that $\la_1 = \frac{1}{3}$, and Lemma \ref{arms} applies.

\end{itemize}

\begin{wrapfigure}[5]{r}{3.1cm}
\vspace{-.25in}
\includegraphics[width=2.8cm]{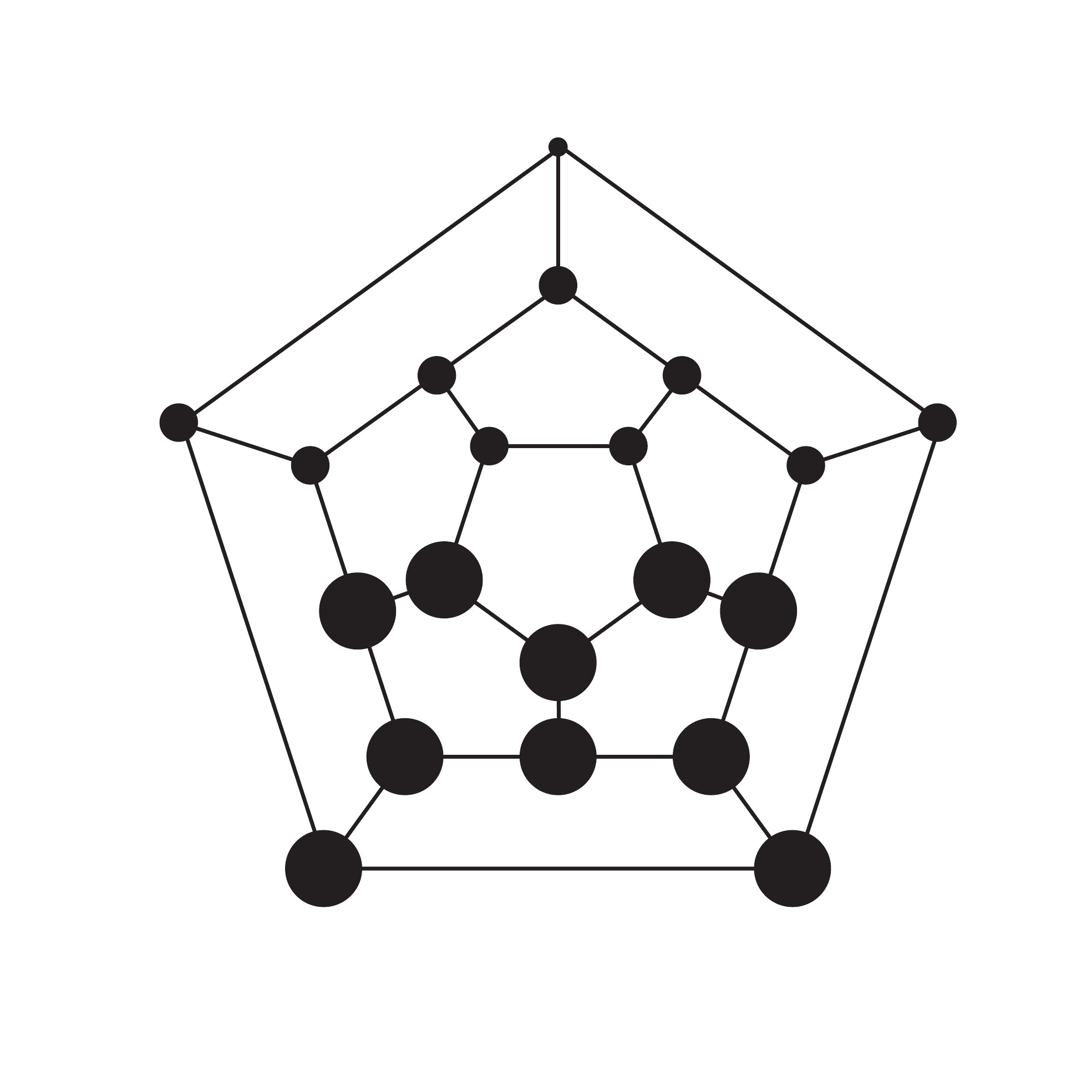}
\end{wrapfigure}

\setlength{\parindent}{13.25pt}

\hspace{0pt} \textbullet \hspace{1.5pt} The dodecahedron, with intersection array $\{3,2,1,1,1;1,1,1,2,3\}$. Lemma \ref{arms} \textcolor{white}{does} \\

\vspace{-28pt}

\hspace{14pt} does not work here, as $\th_1 = \sqrt{5}$ so that $\la_1 = \frac{3-\sqrt{5}}{3} \approx .255$ Using the represent-\\

\vspace{-28pt}

\hspace{14pt} ation to the right, it may be checked that the set indicated by the larger vertices \\

\vspace{-14pt}

\hspace{14pt} gives an upper bound of $h_G \leq \frac{1}{5}$ (in fact, it can be shown that this is the optimal set).

\vspace{0pt}

\begin{itemize}

\item The Desargues graph, with intersection array $\{3,2,2,1,1;1,1,2,2,3\}$. $\th_1 = 2$ so that $\la_1 = \frac{1}{3}$, and Lemma \ref{arms} applies.

\item Tutte's 12-cage, with intersection array $\{3,2,2,2,2,2;1,1,1,1,1,3\}$. Here $\th_1 = \sqrt{6}$ so that $\la_1 = \frac{3-\sqrt{6}}{3} > .183$. Fix any $x \in G$, and let $G'$ be the distance-2 graph induced on $\Ga_6(x)$. Note that the facts $c_6=3$ and $b_5=2$, together with the absence of cycles in $G$ of length 4, imply that $G'$ has valency 3. Since the girth of $G$ is 12, the girth of $G'$ is at least 6, and we may therefore find a path of length 5 in $G'$. Let $S_6$ be this path together with leaves attached to each of the interior points in $G'$, so $S_6$ induces a tree in $G'$ and $|S_6|=8$. Let $S_5 = \Ga_1(S_6) \cap \Ga_5(x)$ and $S_4 = \Ga_1(S_5) \cap \Ga_4(x)$; note that $c_6=3, c_5=1$ and the absence of short cycles in $G$ imply $|S_5|=17$. Furthermore, $c_5=1$ implies $|S_4| =: a \leq 17$. Note that if $a < 17$ then $17-a$ vertices in $S_4$ are adjacent to two vertices in $S_5$, and therefore $E[S_4,\Ga_5(x) \backslash S_5] = 17-2(17-a) = 2a-17$. We now let $S = \{x\} \cup \Ga_1(x) \cup \Ga_2(x) \cup \Ga_3(x) \cup S_4 \cup S_5 \cup S_6$. Since $\{|\Ga_1(x)|, |\Ga_2(x)|, |\Ga_3(x)|\} = \{3,6,12\}$, we have $|S|= 47 + a$. Furthermore, $E[S,S^c] = E[\Ga_3(x),\Ga_4(x) \backslash S_4] + E[S_4,\Ga_5(x) \backslash S_5] + E[S_5,\Ga_6(x)]= (24-a) + (2a-17) + 10 = a + 17$. If $a=17$, then $|S| = 64$ and $\min(|S|,|S^c|) = |S^c| = 62$, and we have $h_G \leq \frac{34}{3 \times 62} <.183 < \la_1$. If $a<17$, then $\min(|S|,|S^c|) = |S| = 47+a$, and so we have $h_G \leq \frac{E[S,S^c]}{3|S|} = \frac{a+17}{3(a+47)}$. This function is increasing in $a$ so we need only check $a=16$, and $\frac{16+17}{3(16+47)} \approx .175<.183$. The result follows.

\item The Biggs-Smith graph, with intersection array $\{3,2,2,2,1,1,1;1,1,1,1,1,1,3\}$. Here $\th_1 = \frac{1+\sqrt{17}}{2}$ so that $\la_1 = \frac{3-\frac{1+\sqrt{17}}{2}}{3} \approx .146$ As is shown in \cite[p. 405]{drgraphs}, this graph can be decomposed into the following equitable partition.

\begin{center}
\includegraphics[height=4.5cm]{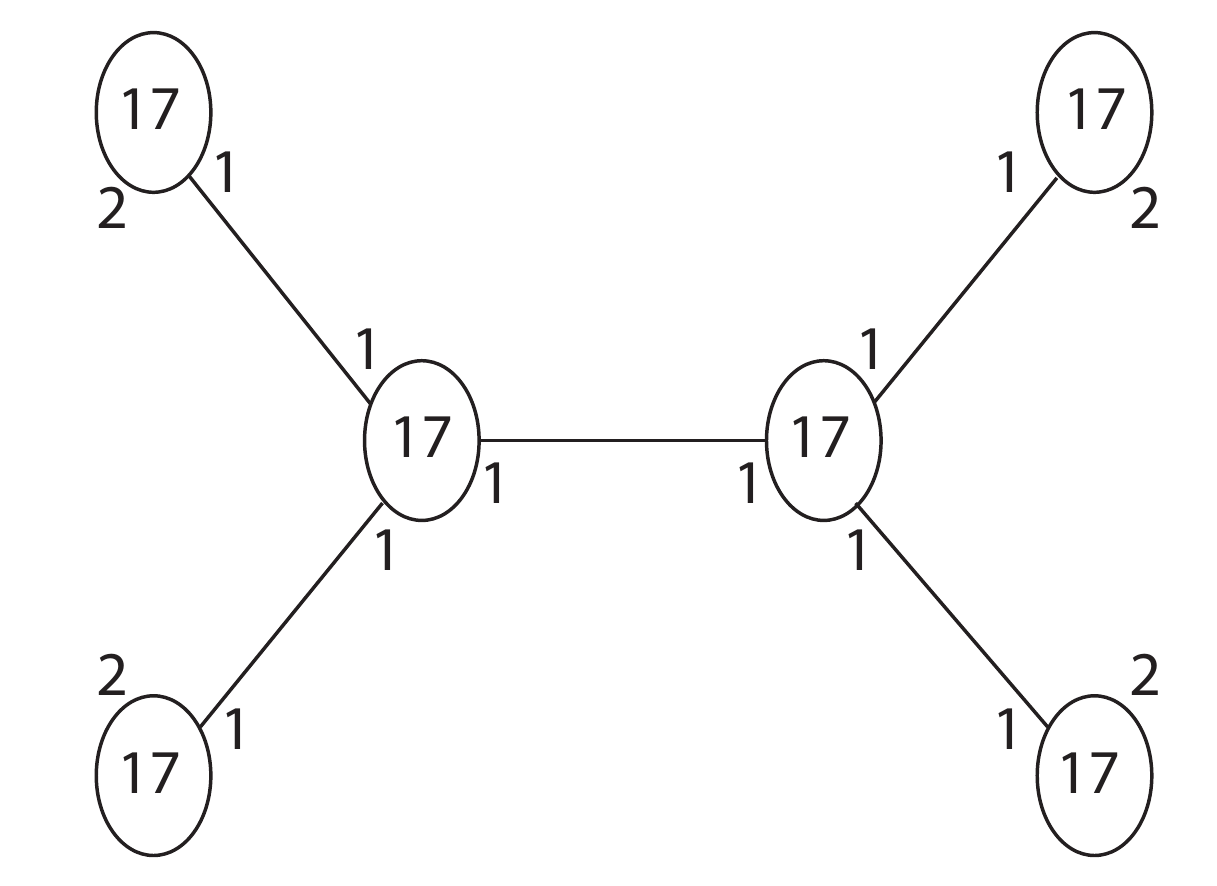}
\end{center}

Take $S$ to be the set of vertices corresponding to the left half of this picture; $S$ has 51 vertices, and $E[S,S^c] = 17$, yielding $h_G \leq \frac{1}{9} < \la_1$. 

\item The Foster graph, with intersection array $\{3,2,2,2,2,1,1,1;1,1,1,1,2,2,2,3\}$. $\th_1 = \sqrt{6}$ so that $\la_1 = \frac{3-\sqrt{6}}{3} \approx .184.$ This graph $G$ has 90 vertices and distribution diagram (\cite[13.2.1]{drgraphs})

\includegraphics[height=2.6cm]{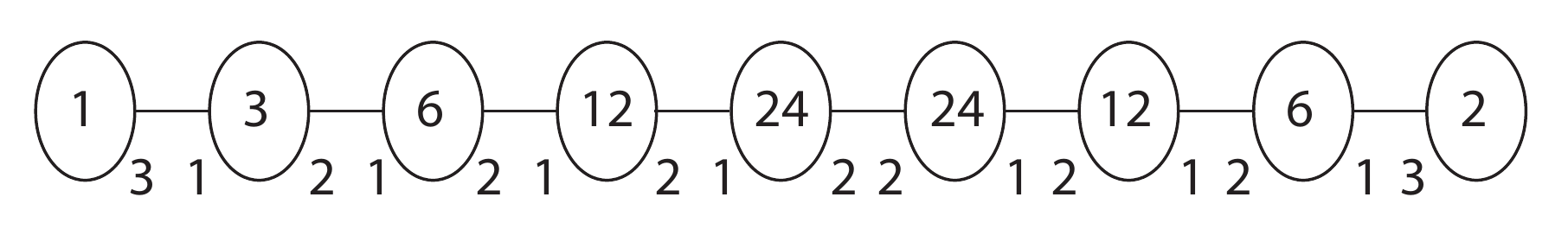}

We are required to find an induced subgraph $G'$ with $|G'|\leq |G|/2$ vertices and $e'\geq \theta_1|G'|/2$ edges, since such a subgraph will satisfy the conditions of Lemma \ref{iran}. We will find a subgraph $G'$ with the following distance distribution diagram for a certain vertex:

\includegraphics[height=2cm]{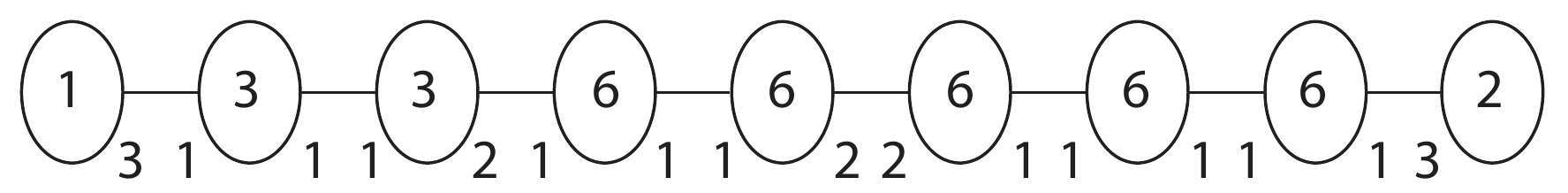}

In this case, $v'=39\leq v/2=45$ and $e'=48\geq \theta_1 |G'|/2\approx 47.77$. 
In order to find such graph, consider the halved Foster graph $F$, which has distance distribution diagram given in \cite[Sec. 13.2B]{drgraphs}. As is shown there, $F$ is an induced subgraph of the Conway-Smith graph, which is the unique connected locally-Petersen graph on 63 vertices. This implies that we can find an induced subgraph of $F$ isomorphic to the Petersen graph. Then the vertex-edge incidence graph of the Petersen graph is an induced subgraph of the Foster graph and has the following distribution distance diagram, where the leftmost balloon contains a vertex corresponding to a vertex in the Petersen graph.

\hspace{3cm} \includegraphics[height=2.4cm]{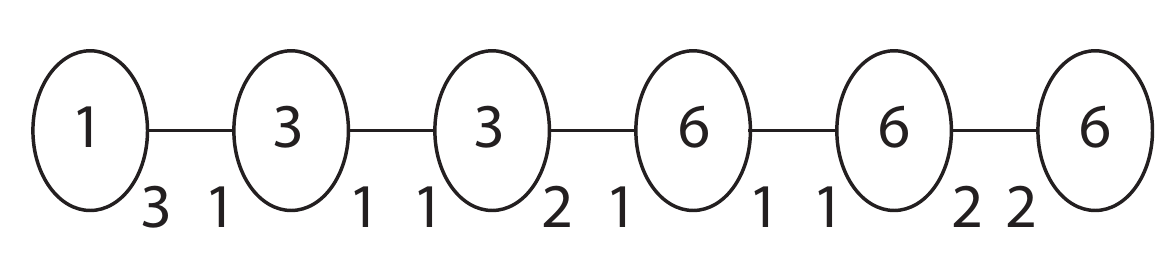}

If consider the two antipodal points to the point in the leftmost balloon, then each has exactly 3 geodesics of length 3 to the rightmost balloon of size 6 (because $c_3=1$), so we may add each of these antipodal points together with each point along these 6 geodesics of length 3. We end up with the required subgraph $G'$. This completes the proof for the Foster graph.

\end{itemize}

As is shown in \cite{broukool} the only distance-regular graphs of valency 4 with diameter greater than 2
are as follows.

\begin{itemize} \label{}

\item $K_{5,5}$ minus a matching, with intersection array $\{4,3,1;1,3,4\}$. $\th_1 = 1$ and thus $\la_1 = \frac{3}{4}$, and Lemma \ref{arms} applies.

\item The nonincidence graph of $PG(2,2)$, with intersection array $\{4,3,2;1,2,4\}$. $\th_1 = \sqrt{2}$ and thus $\la_1 = \frac{4-\sqrt{2}}{4} > \frac{1}{2}$, and Lemma \ref{arms} applies.

\item The line graph of the Petersen graph, with intersection array $\{4,2,1;1,1,4\}$. $\th_1 = 2$ and thus $\la_1 = \frac{1}{2}$, and Lemma \ref{arms} applies.

\item The 4-cube, with intersection array $\{4,3,2,1;1,2,3,4\}$. $\th_1 = 2$ and thus $\la_1 = \frac{1}{2}$, and Lemma \ref{arms} applies.

\item The flag graph of $PG(2,2)$, with intersection array $\{4,2,2;1,1,2\}$. $\th_1 = 1+\sqrt{2}$ and thus $\la_1 = \frac{3-\sqrt{2}}{4} \approx .396$. Since $a_1 = 1$ and $k=4$, for any vertex $x$ we have that $\{x\} \cup \Ga_1(x)$ is two triangles which intersect at a point. We may therefore take $S$ to be a set of $3$ intersecting triangles, which will have $7$ vertices and $E[S,S^c]=10$, giving $h_G \leq \frac{10}{4 \times 7} \approx .357.$

\item The incidence graph of $PG(2,3)$, with intersection array $\{4,3,3;1,1,4\}$. $\th_1 = \sqrt{3}$ and thus $\la_1 = \frac{4-\sqrt{3}}{4} > \frac{1}{2}$, and Lemma \ref{arms} applies.

\item The incidence graph of $AG(2,4)$ minus a parallel class, with intersection array $\{4,3,3,1;1,1,3,4\}$. $\th_1 = 2$ and thus $\la_1 = \frac{1}{2}$, and Lemma \ref{arms} applies.

\item The odd graph $O_4$, with intersection array $\{4,3,3;1,1,2\}$. This was shown to be OK in Section \ref{infam}; alternatively $\th_1 = 2$ and thus $\la_1 = \frac{1}{2}$, and Lemma \ref{arms} applies.

\item The flag graph of $GQ(2,2)$, with intersection array $\{4,2,2,2;1,1,1,2\}$. $\th_1 = 3$ and thus $\la_1 = \frac{1}{4}$. Since $a_1 = 1$ and $k=4$, for any vertex $x$ we have that $\{x\} \cup \Ga_1(x)$ is two triangles which intersect at a point, and it also follows that any edge is contained in a unique triangle. Take two points at distance 4 in the graph, and then consider two disjoint 4-paths between the points. Let $S$ be the set of all points contained in these paths, together with their neighbors. Since each edge is contained in a triangle, $S$ will induce an octagon with each edge replaced by a triangle. Thus, $|S|=16$, and $E[S,S^c] = 16$. We see $h_G \leq \frac{16}{4 \times 16} = \frac{1}{4}$.

\item The doubled odd graph $(dO)_4$, with intersection array $\{4,3,3,2,2,1,1;1,1,2,2,3,3,4\}$. This was shown to be OK in Section \ref{infam}. 

\item The incidence graph of $GQ(3,3)$, with intersection array $\{4,3,3,3;1,1,1,4\}$. $\th_1 = \sqrt{6}$ and thus $\la_1 = \frac{4-\sqrt{6}}{4} \approx .388$. This graph is bipartite, and has distribution diagram

\hspace{3cm} \includegraphics[height=2.4cm]{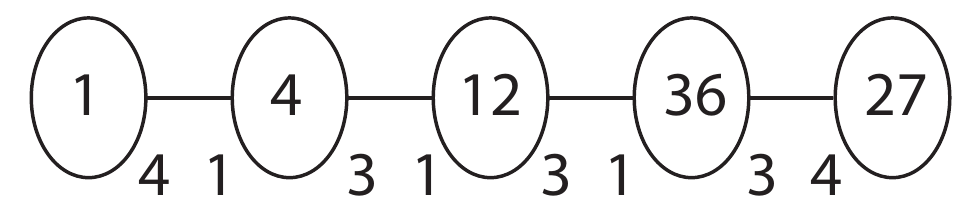}

For fixed $x \in G$ we may choose three points $y_i$ for $i=1,2,3$ in $\Ga_4(x)$ which are mutually of distance 4 from each other; this is possible, since $|\Ga_2(y_i)\cap \Ga_4(x)|=8$. We then let $S=\{x\}\cup \Ga(x)\cup \Ga_2(x) \cup \Big(\cup_{i=1}^3 (\{y_i\}\cup \Ga(y_i))\Big)$. It can be then be checked that $|S|= 32 \leq \frac{|G|}{2}$ and $E[S,S^c]=48$, so that $h_G \leq \frac{48}{4 \times 32} = .375$. We note that while this graph is known to exist, it is not yet known whether it is unique; however, our proof uses only properties of the intersection array, and therefore applies to all such graphs, as their intersection arrays coincide.



\item The flag graph of $GH(2,2)$, with intersection array $\{4,2,2,2,2,2;1,1,1,1,1,2\}$. $\th_1 = 1+\sqrt{6}$ and thus $\la_1 = \frac{3-\sqrt{6}}{4}$. We have not been able to determine whether or not this graph is OK.

\item The incidence graph of $GH(3,3)$, with intersection array $\{4,3,3,3,3,3;1,1,1,1,1,4\}$. $\th_1 = 3$ and thus $\la_1 = \frac{1}{4}$. We have not been able to determine whether or not this graph is OK.

\end{itemize}

\section*{Acknowledgements}


We would like to thank the anonymous referees for their careful reading. J.H.K. has been partially supported by the National Natural Science Foundation of China (Grants No. 11471009 and No. 11671376) and by the Anhui Initiative in Quantum Information Technologies (Grant No. AHY150200). G.M. has been partially supported by the Australian Research Council (Grants DP0988483 and DE140101201). Z.Q. has been partially supported by the National Natural Science Foundation of China (Grant No. 11801388).

\bibliographystyle{alpha}
\bibliography{drgraphs}

\end{document}